\documentclass[11pt]{article}
\usepackage{amsmath,amsthm,amssymb}
\usepackage{graphicx, color}
\usepackage{psfrag}
\usepackage{epsf}
\usepackage{array}
\usepackage[title]{appendix}
\usepackage[mathscr]{eucal}
\usepackage[T1]{fontenc}
\usepackage[utf8]{inputenc}
\usepackage[flushleft]{threeparttable}
\usepackage{float}
\usepackage{booktabs,caption,fixltx2e}
\usepackage{authblk}

\usepackage{booktabs}
\usepackage{caption}
\usepackage{float}
\usepackage{subcaption}

\usepackage{amsmath,amsthm,amssymb}
\usepackage{graphicx, color}
\usepackage{psfrag}
\usepackage{epsf}
\usepackage{array}
\usepackage[title]{appendix}
\usepackage[mathscr]{eucal}
\usepackage[T1]{fontenc}
\usepackage[utf8]{inputenc}
\usepackage[flushleft]{threeparttable}
\usepackage{float}
\usepackage{booktabs,caption,fixltx2e}
\usepackage{authblk}

\usepackage{booktabs}
\usepackage{caption}
\usepackage{float}
\usepackage{subcaption}

\newtheorem{theorem}{Theorem}

\newtheorem{criterion}[theorem]{Criterion}

\newtheorem{remark}{Remark}

\newtheorem{assumption}{Assumption}

\begin{document}
	
	
	
	%
	%
	
	\pagestyle{plain}
		
		\title{A Decomposition Property for an  $M^{X}/G/1$ Queue with Vacations }
		\author{Igor Kleiner\thanks{Department of Statistics, Haifa University, Haifa, Israel (igkleiner@gmail.com)},
			\hspace{2mm}      
			Esther Frostig\thanks{Department of Statistics, Haifa University, Haifa, Israel (frostig@stat.haifa.ac.il)
			\\
		The research of Esther Frostig is partly funded by ISF (Israel Science Foundation), Grant
		Grant  1999/18).}
		\hspace{2mm}and David Perry
	\thanks{ Holon Institute of Technology, P.O. Box 305, Holon 5810201, Israel (davidper@hit.ac.il)
\\
The research of David Perry is partly funded by ISF (Israel Science Foundation), Grant 3274/19.}
}
		
		%
		%
		
		%
	

		\maketitle              
	
		\begin{abstract}

	 We introduce a  queueing  system that alternates between two modes, so- called  {\it working mode} and {\it vacation mode}. During the working mode the system runs as an $M^{X}/G/1$ queue. Once the number of customers in the working mode drops to zero the vacation mode begins. 
			During the  vacation mode the system runs as a general queueing   system (a service might be included) which is different from the one in the working mode. 
			The vacation period ends in accordance with  a given stopping rule, and then a random number of customers are transferred to the working mode. 
			For this model we show  that the conditional probability 
			generating function  of the number of customers given that the system is in the working mode is a product of three terms.
			 This decomposition result puts under the same umbrella some  models that have already been introduced in the past as well as  some new models. 
			%
			%
			
	\vspace{2mm}
		\end{abstract}
		
		\section{Introduction}

	Consider a queueing system alternating between two modes: 
	{\it working mode} and {\it vacation mode}.
	During the working period the system runs as a stable $M^{X}/G/1$ queue. The working period   terminates   whenever the number of customers drops to $0$. Then, the system starts the vacation mode and it stays there until the next  transition back to the working mode. This  occurs according  to a certain stopping rule.
	
	In this study, the code name {\it vacation} 
	is interpreted in the broad sense; namely, during the  vacation period the system  runs as a general queueing system, possibly including general services, batch arrivals, arbitrary number of servers, abandonment,  disasters,  etc.

	It is assumed that the vacation is  terminated with at 
	least one customer, which means that the emptiness period is included in the vacation period.   At the transition  moments   all the customers,   including those being  served,  are transferred  altogether   to the working mode, such that the  service times of the customers being served  are lost;  their  
	services  start from scratch   in  the working mode,  possibly with different service regime.


	The  working period  starts with  a random number of customers (those who have just been transferred from the vacation mode) and ends whenever the number of customers  drops to $0$. Thus, during the working mode the system runs as a modified  $M^{X}/G/1$ queueing system, where the batches arriving  during the busy period are i.i.d., but the first batch is of  a different distribution from the others. The latter distribution  is that of  the number of customers that are transferred  from the vacation mode to the working mode.

	

	
	The decomposition  property studied here  says that the conditional steady state probability generating function (PGF)
	of the    number of customers, given the  working mode,   is a product of three terms:  the conditional PGF of the number of customers in a regular $M^{X}/G/1$ queue, given that the server is busy,  the equilibrium PGF of the number of customers transferred from the vacation mode to the working mode and    the reciprocal of the  equilibrium PGF of the number of customers in  a batch.   
	
	This study generalizes the decomposition property by considering an $M^{X}/G/1$ queueing system  in the working mode and a general queueing system in the vacation mode.
	%
	
	The decomposition property   has been studied over the last three decades as a  useful tool for analyzing  specific  queueing systems with vacations.  Two vacation variants  are  considered: (i) an {\it idle vacation} without services  and (ii) a {\it working vacation} with services. 	
	 
	The  decomposition result introduced here  characterizes   a general framework that enables to analyze many  vacation queues that are studied  in the literature as special cases. Also, it enables to analyze more general models. 
	A comprehensive review of  this property for queues with vacation until 1986  can be found in  Doshi \cite{Doshi} and the references therein.

Usually, in the models of type (i) the system runs as an $M/G/1$ queue during the working period.  Altman and Yechiali \cite{Altman} studied a model with exponentially distributed vacation period and with reneging during vacations. Adan et al. \cite{IvoNew} study models with batch reneging (synchronized reneging)  during vacation, where  the vacation period is generally distributed. 
Cooper  and Fuhrmann \cite{CooperFDec} prove the   decomposition property for a class of  $M/G/1$ type  queues. 
Shanthikumar \cite{Shanti}  proves  the decomposition result for the number of customers in the system at departure times in  the working mode under more general assumptions than in \cite{CooperFDec}.
             Baba \cite{Baba87}  and Choudhury \cite{Chou} consider the distribution of the number of customers in an  $M^X/G/1$  vacation queue with general vacation times. Baba \cite{Baba87} uses the supplementary variables technique and \cite{Chou} demonstrates the   decomposition property.                         
             Mytalas and Zazanis \cite{Mytalas} study an $M^{X}/G/1 $ type  system with disasters with two types of vacations: repair after disasters, and multiple vacations after reaching state 0 due to a service completion.  They  use the  supplementary variables   methodology to obtain  explicit  expressions for the partial PGFs of the  number of customers in the working period, the  number of customers when the system is being  repaired, and the  number of customers during the  vacation. They also find  the  PGF of the number of customers at  departure epochs.
		      Levy and Yechali \cite{Yech75} analyze the $M/M/s$ vacation queues. The latter two studies  do not prove or use any decomposition result. 

		 In type (ii) models 		 
		  the vacation period  terminates with at least one customer. The service of the customers being served at the transition epoch from the vacation mode to the working mode is  lost and their  service starts from scratch in the working mode. 
		   Liu et al. \cite{Liu} introduce  a decomposition result when the system runs as an M/M/1 queue both in the working and in  the vacation modes, but with different   service rates.  The duration of a vacation is exponential.	   
		  Li  et al. \cite{Li} consider an M/G/1 queue, both in the working and in the  vacation modes. However, the service distribution in the vacation mode is different from that in the working mode. The vacation period is  
		    exponentially
		   distributed. 
		    They use the  matrix-analytics method and Shanthikumar \cite{Shanti} to obtain  the stationary distribution of the number of customers at  departure epochs.  
		   Gao and Yao \cite{Gao} assume that both in the working and in the  vacation modes the system runs as an $M^{X}/G/1$ queue but with different service distributions. The total vacation period is as follows. First, an exponentially distributed random variable is sampled.   If it terminates when the system is not empty, there is a transition to the working mode. Otherwise,  with probability $p$  the vacation terminates with zero  customers and the working mode starts with zero customers, alternatively,  with probability $1-p$ another exponential  random variable is sampled and so forth.   The number of returns to state 0 in the vacation is bounded.  Note, that in their model the system can be empty  in the vacation mode  and   in the working mode. 		   
		    By using the supplementary variables technique   they find the   distribution  of the number of customers in the system  in steady state as well as additional quantities of interest. 

	 In section 2, we introduce the model and prove the main decomposition result.  In section 3 we apply the 
		 decomposition result to analyze an $M^{X}/G/1 $ queue in the working mode and general types of queueing  models that run during the vacation mode.

	\section{The Decomposition Result}
We start with  a detailed description of the model  and introduce some assumptions and notations.  Then  the  decomposition result is proven. 
	\subsection{Model Description}\label{secdescription}
	 
 During the working mode the system behaves as an
$M^{X}/G/1$ queueing system, and during the vacation mode  the system's  behavior  is arbitrary  general. Throughout, we assume that the system is in steady-state.
	
	 Upon  termination of the vacation mode  the service time of the customers in service is lost. These customers restart their service (probably with different service distribution) at the working mode. 
 Denote by $Y$ the generic random variable  that describes   the number of customers   transferred to the working mode at the end of the vacation mode and let  
$\psi_i=P(Y=i),i=1,2,3,...$ be its probability function. Denote by $\Psi(z)$ its PGF

\begin{equation}
\Psi(z)=\sum_{i=1}^{\infty}\psi_{i} z^{i}.
\label{psynew}
\end{equation}	

The vacation terminates with at least one customer. 
Note that $\Psi(z)$  is not yet given in terms of known functionals; obviously, it is determined  in accordance with the specific model that is run in the vacation mode.

In the  working mode the system's behavior is  that of a modified   $M^X/G/1$ queueing system, since the distribution of the number $Y$  of  transferred customers, is   different
from the  given  distribution of the regular batch sizes.    During the busy (working) period of the latter $M^X/G/1$ queueing system batches arrive according to a $Poisson$ process with rate $\lambda$. 	Let $S$ be the generic service time of a customer whose distribution is $F_S$ and Laplace Stieltjes transform (LST) $\widetilde
{F}_S(\alpha)$.	
Let $B$ be  the generic batch size. Let    $b_i=P(B=i),i=1,2,...$ be its probability vector and 
\begin{equation}
	B(z)=\sum_{i=1}^{\infty}b_i z^i
\end{equation}
be its PGF.

	Once the number of customers in the working mode drops to $0$ the vacation mode starts. 
	
	To summarize the above, we indicate that in the	working mode the system evolves as a modified 
 $M^{X} /G/1$ system, where the busy period starts with  $Y$ customers and all the  batches arriving during the busy period are i.i.d. distributed as $B$.  
Denote this system by $[M^{X}/G/1]^{Y}$.	
 Specifically, it is assumed that
		 the expected vacation period is 
		finite, $E[Y]<\infty$, 
		the queueing system in vacation mode is in steady state and 
	 $\lambda E[S]E[B]<1$. 
	So that the system is regenerative. 
	Let $T$ be the length of a cycle, that is, the duration of the vacation mode and the working mode that follows. 
	
We say that the system is in state  $(0,j)$ , $j=0,1,2,...$  if there are $j$ customers in the vacation mode. Let $p_{(0,j)},j=0,1,2,...$ be the steady--state probabilities that there are $j$ customers in the vacation mode and let  
 \begin{equation}
G_0(z)=\sum_{i=0}^{\infty}p_{(0,i)} z^i
\end{equation}
be its PGF.  
Similarly, we say that the system is in state $(1,j)$,  if there are $j$ customers in the working mode.	Let $p_{(1,j)},j=1,2,...$, $j=1,2,...$ be the steady--state probabilities that the system is in state $(1,j)$ and let 
	 \begin{equation}
	 G_1(z)=\sum_{i=1}^{\infty}p_{(1,i)} z^i
	 \end{equation}
	 be its PGF.
	 Then
	  \[p_{0.}=\sum_{i=0}^{\infty}p_{(0,i)}\]
is the steady--state probability that the system is in the vacation mode and  \[p_{1.}=1-p_{0.}=\sum_{i=1}^{\infty}p_{(1,i)}\]
 is the steady--state probability that the system is in the working mode. 
	 Finally, let $\widetilde{G}_{1}(z)$ be the conditional PGF of the number of customers, given that the system is in the working mode. That is 
	 \[\widetilde{G}_{1}(z)=\frac{ G_1(z)}{p_1.}.\]
	 
Clearly,	 $\widetilde{G}_{1}(z)$ is  the PGF of the number of customers, given that the $[M^{X}/G/1]^{Y}$ system is not idle.
	
For a non-negative integer random variable $U$ with finite expectation and  PGF  $U(z)$, we denote by $U^{e}$ an integer valued random variable with probability function
\begin{equation}\label{eq.equi}
P(U^{e}=i)=\frac{P(U>i)}{E[U]}.
\end{equation}
That is, $U^{e}$ is   the equilibrium random variable associated with  $U$; its PGF $U^{e}(z)$  is
\begin{equation}\label{equilibriumPGF}
 U^{e}(z)=\frac{1-U(z)}{E[U](1-z)}.
 \end{equation}

\subsection{Decomposition }
Consider a regular   $M^X/G/1$ queue where all the batches arriving at the system are distributed as $B$ and let 
 $P(z)$ be the PGF of the number of customers in steady--state.   By Tijms \cite{Tijms} 
\begin{equation}
P(z)=(1-\lambda E[S]E[B])\frac{1-\lambda \alpha(z)(1-B(z))}{1-\lambda\alpha(z)\frac{1-B(z)}{1-z}},
\end{equation}
where 
\begin{equation}
\alpha(z)=\int_0^{\infty} e^{-\lambda (1-B(z))t}(1-F_S(t))dt.
\end{equation}
Clearly,   $\rho=\lambda E[S]E[B]<1$ and it is the probability that the server is busy. 
Let $\widetilde{P}_{1}(z)$ be the conditional PGF of the number of customers in the $M^X/G/1$ queue    given that the server is busy. Then
\begin{eqnarray}
 \widetilde{P}_{1}(z)&&=\frac{P(z)-(1-\rho)}{\rho}\nonumber\\ 
&&= \frac{1-\rho}{\rho} \frac{\lambda z \alpha(z)\frac{1-B(z)}{1-z}}{1-\lambda \alpha(z)\frac{1-B(z)}{1-z}}\nonumber\\.\label{tildeP1}
&&\nonumber\\
&&=\frac{1-\rho}{\lambda E[S]}\frac{\lambda z \alpha(z)B^{e}(z)}{1-\lambda \alpha(z)\frac{1-B(z)}{1-z}}
 \end{eqnarray} 


We are now in a position  to introduce  the decomposition result. To this end   we introduce (similarly to 9.3.1 in \cite{Tijms})
 $a_{j}$--the expected amount of  time that  $k+j$ customers are present  during a given  service time that starts with $k$ customers. This quantity is independent of $k$ so that
\begin{equation}\label{eq.an}
	a_j=\int _0 ^\infty r_j(t)(1-F_{S}(t))dt,\,\,\,j=0,1,...
\end{equation}
where
\begin{equation}
	r_{j}(t)=P(j \text{ customers arrive in }(0,t) \text{ in the working mode}).
\end{equation}

\begin{theorem}
\label{theorem.main}
	\begin{eqnarray}\label{decomposition}
\widetilde{G}_{1}(z)=\widetilde{P}_{1}(z)(B^e(z))^{-1}\bullet\Psi ^e(z),   
\end{eqnarray}
where $B^e(z)$ and $\Psi ^e(z)$ 
are the PGFs of the equilibrium distribution functions of $B$ and $Y$, respectively. 
\end{theorem}
		The reason for  the appearance of  $(B^{e}(z))^{-1}$  on the right hand side of (\ref{decomposition}) can be motivated as follows.
		  The numerator  of (\ref{tildeP1}) contains  $B^{e}(z)$ which is the PGF of the equilibrium  batch size.
		In the regular $M^{X}/G/1$ queue, all the   batches are i.i.d.,  while  in our  $M^{X}/G/1$  modified model (with vacations)   the PGF of the batch that starts the working mode is   $\Psi(z)$. As a result, $B^{e}(z)$ is replaced by $\Psi^{e}(z)$ in (\ref{decomposition}) .
\begin{proof}
  We start with  a recursive scheme to obtain $p_{(1,j)}$, $j=1,2,...$. 
\begin{equation}\label{recrul}
p_{(1,j)}=\frac{\sum _{s=1}^{j} \psi _s a_{j-s}}{ET} +  \sum_{k=1}^{j}a_{j-k} \left(\frac{ P(Y>k)}{E[T]}+\lambda \sum _{i=1} ^k p_{(1,i)} P(B>k-i)  \right).
 \end{equation}
To explain the recursive scheme  (\ref{recrul}),
recall that $T$  is  the length of the cycle, i.e. the length of the vacation period plus the  working period.  Denote by $T_j$ the duration  of time in one cycle in which there are j customers in the working mode. Then by renewal theory
\begin{equation}
p_{(1,j)}=\frac{E[T_j]}{E[T]} \label{p1j}.
\end{equation} 

Let $N_k$ be the number of service completions in the  working mode in one 
cycle  in which $k$ customers are left behind. In other words, 
$N_k$ is the number of transitions from state $(1,k+1)$ to state $(1,k)$ in one cycle.  
 Recall that for $j\geq k$, $a_{j-k}$ is the expected time during a service time that starts with $k$ customers during which  $j$ customers are present.
By Wald's identity
\begin{equation}\label{ETJ}
E[T_j]=\sum_{s=1}^j \psi_s a_{j-s} + \sum_{k=1}^j E[N_k] a_{j-k}. 
\end{equation}

To obtain $E[N_k]$ we apply the level-crossing argument: $N_{k}$,
the number of down-crossings from state $(1,k+1)$ to state $(1,k)$,
   is equal to the number of batch arrivals in one cycle  that bring  the system to state $(1,h), h>k$ (the number of up-crossings of level $k$) at arrival epochs. This event occurs either when the number of customers transferred from the vacation mode is bigger than $k$, or when the size of a batch that arrives during $T_{j}$ is bigger than $k-j$,  $j=1,...,k$. As a result

\begin{equation}
E[N_k] = P(Y>k) + \lambda \sum_{j=1}^k E[T_j] P(B>k-j), \,\,k=1,2,.... \label{ENK}
\end{equation} 
Now, substituting (\ref{ENK}) in (\ref{ETJ}), dividing  by $E[T]$  and applying  (\ref{p1j}) we obtain  (\ref{recrul}).

Finally, the decomposition result (\ref{decomposition}) is obtained by multiplying both sides of  (\ref{recrul}) by $z^j$, summing for $j=1,2,...$ and dividing by  $p_{1.}$. 
Obviously,  the expected time of the working mode in one cycle is
\[\frac{E[S]E[Y]}{1-\lambda E[B]E[S]},\]
and by renewal theory 
\[p_{1.}=\frac{1}{E[T]}\frac{E[S]E[Y]}{1-\lambda E[B]E[S]}.\]
%
\end{proof}

	
{\large{\bf Special cases:} }
\begin{enumerate}
	\item $Y$ is distributed as $B$. This case occurs  when the vacation period coincides with the idle period. Then  $B^{e}(z)=\Psi^{e}(z)$
	 and $\widetilde{G}_{1}(z)=\widetilde{P}_{1}(z)$ is given in (\ref{tildeP1}).
	
	\item  $B=1$. In this case $G_{B}(z)=1$ and $E[B]=1$.
	Then (\ref{decomposition}) is given by
	\begin{equation}\label{eq.decomM/G/1vac1}
	\widetilde{G}_{1}(z)=\underbrace{\frac{1-\lambda E(S)}{\lambda E(S)}\frac{\lambda \alpha(z) z }{1-\lambda \alpha(z)}}_{I}\underbrace{\frac{(1-\Psi(z))}{E(Y)(1-z)}}_{II},
	\end{equation}
where 
\begin{equation}
\alpha(z)=\int_0^\infty (1-F_{S}(t)) e^{-\lambda(1-z)t}dt=\frac{1-\tilde{F}_S(\lambda(1-z))}{\lambda(1-z)}.
\end{equation}
	Note that the right hand side  of (\ref{eq.decomM/G/1vac1}) is a product of two terms -- I and II. Term I is the conditional PGF of the number of customers in a regular M/G/1 queue   given that the server is busy. Term II is the PGF of $Y^{e}$. 
	This  result is introduced  also in  Equation 5.4 of Adan et al. \cite{IvoNew}.
	
	\item  $Y=1$. This special case is artificial, but it is important for future use. In this case,  each working mode starts with one customer but during the working mode  customers arrive in batches. By (\ref{equilibriumPGF}) $\Psi^{e}(z)=1$.  By substituting (\ref{tildeP1})  in  (\ref{decomposition}) we obtain
	\begin{eqnarray}
\widetilde{G}_{1}(z)	&&=\widetilde{P_1}(z)(B^e(z))^{-1}\nonumber\\
	&&\nonumber\\
	&&=\frac{1-\lambda E(B)E(S)}{\lambda E(S)}\frac{\lambda \alpha(z) z }{1-\lambda \alpha(z)\frac{1-B(z)}{1-z}}.\label{eq.mxg11}	
	\end{eqnarray}
\end{enumerate}

\begin{remark}
The term on the left of the $'\bullet'$ in (\ref{decomposition}) is the same as 
that of (\ref{eq.mxg11}) and the term on the 
right of the $'\bullet'$ is the PGF of the equilibrium distribution of the number of customers at the end of vacation. It follows from (\ref{decomposition}) that the conditional distribution of the  number of customers given that the system is in the working mode  is a convolution of the conditional distribution  of the number of customers given that the working mode starts  with one customer   and  the equilibrium distribution of the number of customers at the end of vacation. The latter statement is  an  additional interpretation of the decomposition result. 
\end{remark}
 
 \begin{remark}\label{rem.2}
 	 	An alternative proof of Theorem \ref{theorem.main} is based on
 	 	 sample path analysis in one cycle. Consider the sample path of the number of customers in the working mode. That is, from the sample path of the original process the vacation periods are deleted and the working   periods are glued together. All the working  periods are i.i.d., each of them starts with a random  number of customers, distributed as $Y$. The length of the working period $\cal{T}$ is presented as the sum of $Y$ sub-busy periods
 	\begin{equation}
 	\mathcal{T}=C_{1}+\cdots+ C_{Y},
 	\end{equation}
 	where the sub-busy period $C_{j}$ is the time it takes to drop from $j$ customers to $j-1$ customers, $j=1,2,..,Y$.  $C_{1},C_{2},...$ are i.i.d. and independent of $Y$, see Figure (\ref{figst1}).   $C_{j}$ is distributed as the busy period in an $[M^{X}/G/1]^{1}$ system. Let $L$ be the generic number of customers in an $[M^{X}/G/1]^{1}$ system   given that the server is busy. The PGF of $L$ is  given in (\ref{eq.mxg11}). 
  	  Let $W^{Y}$ be the number of customers in the system at the beginning of such a sub-busy period. 	Since  the sub-busy periods are glued together it follows by PASTA that the  number of customers seen at the arrival  is $W^{Y}-1+L$, where $W^{Y}$ and $L$ are independent. Furthermore, by  the {\tt inspection paradox} (Ross \cite{Ross2} p.117) the PGF of $W^{Y}-1$ is $\frac{1-\Psi(z)}{E(Y)(1-z)}$, and  (\ref{decomposition}) obtained by using (\ref{eq.mxg11}).
 	
\begin{figure}[h!]
			\centering
			\includegraphics[width=1\textwidth]{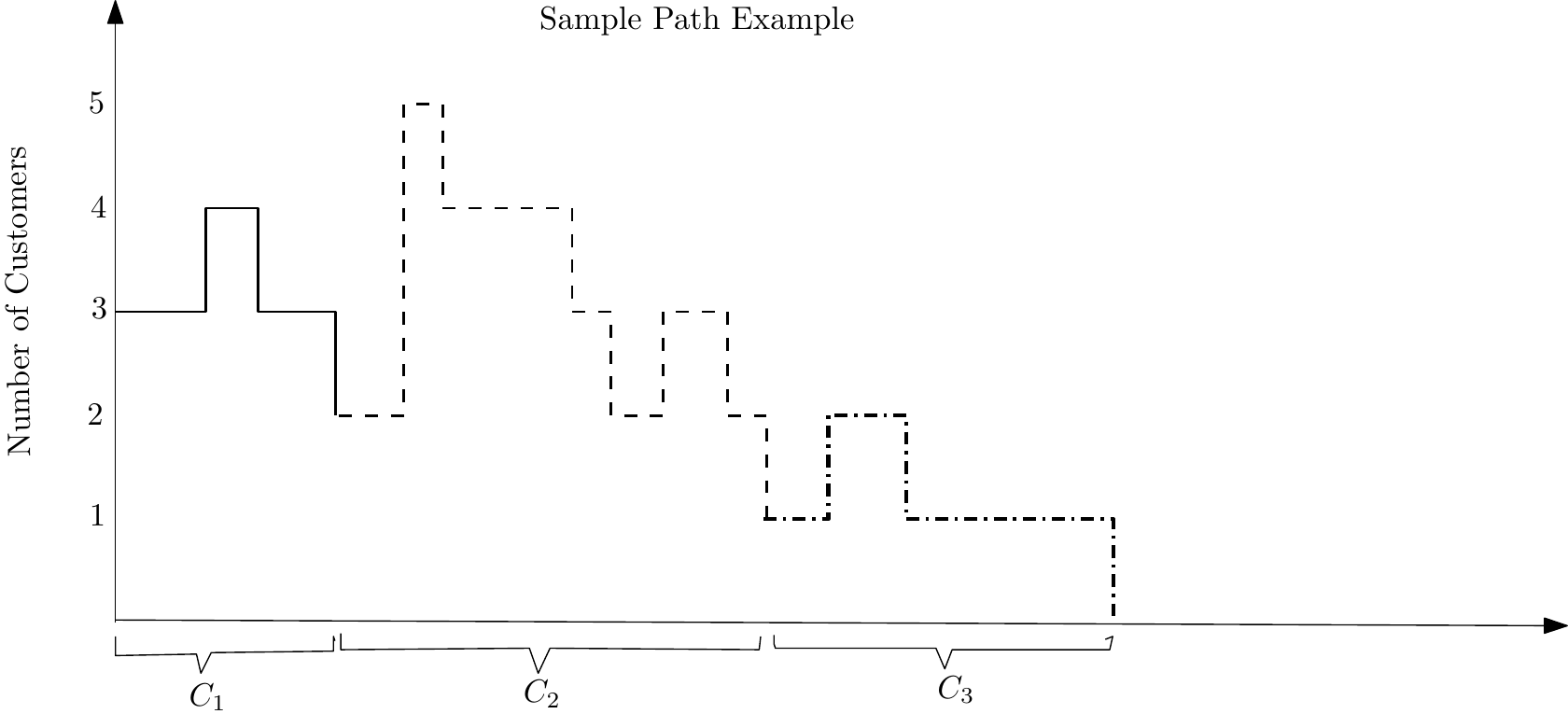}
			\caption{}
			\label{figst1}
\end{figure}
 \end{remark}
		
Motivated by  Remark \ref{rem.2} one can think of   a certain variation of a {\it quasi} $M^{X}/G/c$ model with vacations  for which  a similar decomposition property holds. In this variant  there are $c$ servers in the working mode. Each busy period in the working mode starts with a random number of customers $Y$ transferred from the vacation mode. We refer to  these customers as  {\it prime customers}. The first prime  customer starts a busy period of an  $M^{X}/G/c$ model while all the rest of the $Y-1$   prime customers (if any) wait in line until the service  completion of the prime customer and all his descendants being served. 
Then the next prime customer (if any) starts a stochastically equal  busy period and so forth until the number of customers in the working mode drops to zero.    The latter moment is the beginning of   the  vacation mode of a new cycle.  
Denote this system by  $[\widetilde{M^{X}/G/c}]^{Y}$ and let $G_{1}^{[\widetilde{M^{X}G/c}]^{Y}}$ be the PGF of  the number of customers given that the system is in the working mode. Then, similarly to  (\ref{theorem.main}) we have 
 
 \begin{equation}
 G_{1}^{[\widetilde{M^{X}G/c}]^{Y}}(z)=G_{1}^{[\widetilde{M^{X}G/c}]^{1}}(z)
 \bullet \frac{1-\Psi(z)}{E(Y)(1-z)},
 \end{equation}
which means that the conditional  number of customers given that  the system is in the  working mode   with general $Y$ has the same distribution as that of the convolution of the   conditional  distribution of the number of customers given that the system is in the working mode  with $Y=1$ and the equilibrium distribution of $Y$. 

\section{ Applications}
Throughout this section we assume that during the working mode the system's behavior is that of  an $M^{X}/G/1$ queue as described in Section \ref{secdescription}. We introduce  several examples  with different  behaviors in the vacation mode, for which we obtain the PGF $\Psi(z)$ of the number of customers  transferred from the vacation mode to the working mode. According to Theorem \ref{theorem.main} the PGF of the number of customers given that the system is in the working mode is obtained.
\subsection{Multiple vacations}\label{sec.multiplevacation}
Once the number of customers in the working mode drops to zero
the server takes a vacation for a random time  $V_{1}$, whose distribution is  $F_{V}$ with LST  $\widetilde{F}_{V}(s)$. If at time $V_{1}$ there is at least one customer  the vacation terminates and the customers accumulated during $V_{1}$ are transferred to the working mode. Otherwise, the server starts another vacation  $V_{2}$ independent of the other vacations distributed as $V$, and so forth. The number of customers at the end of vacation mode -- $Y$ is distributed as the number of customers arriving during $V$ given that at least one customer arrived during this period. During the vacation the customers arrive in batches   according to a Poisson process with rate $\lambda_{v}$ and the batch sizes are i.i.d  distributed as $B_{v}$  with PGF $G_{B_{v}}(z)$. In this case the PGF $\Psi(z)$ of $Y$, is
\begin{equation}\label{psimultiple}
\Psi(z)=\frac{\tilde{F}_{V}(\lambda_{v}(1-G_{B_{v}}(z)))-\tilde{F}_{v}(\lambda_{v})}{1-\tilde{F}_{V}(\lambda_{v})}.
\end{equation}
 The PGF of the number of customers given that the system is in the   working mode is obtained by substituting (\ref{psimultiple}) in (\ref{theorem.main}). This result is  obtained  by Lee et al. \cite{LeeChae}, for the case when the vacation terminates with at least one customer ($N=1$ policy).

\subsection{A Markovian  vacation mode   }
\subsubsection{A Markovian vacation with batch arrivals and balking}
\begin{figure}[h!]
	\centering
	\includegraphics[width=1\textwidth]{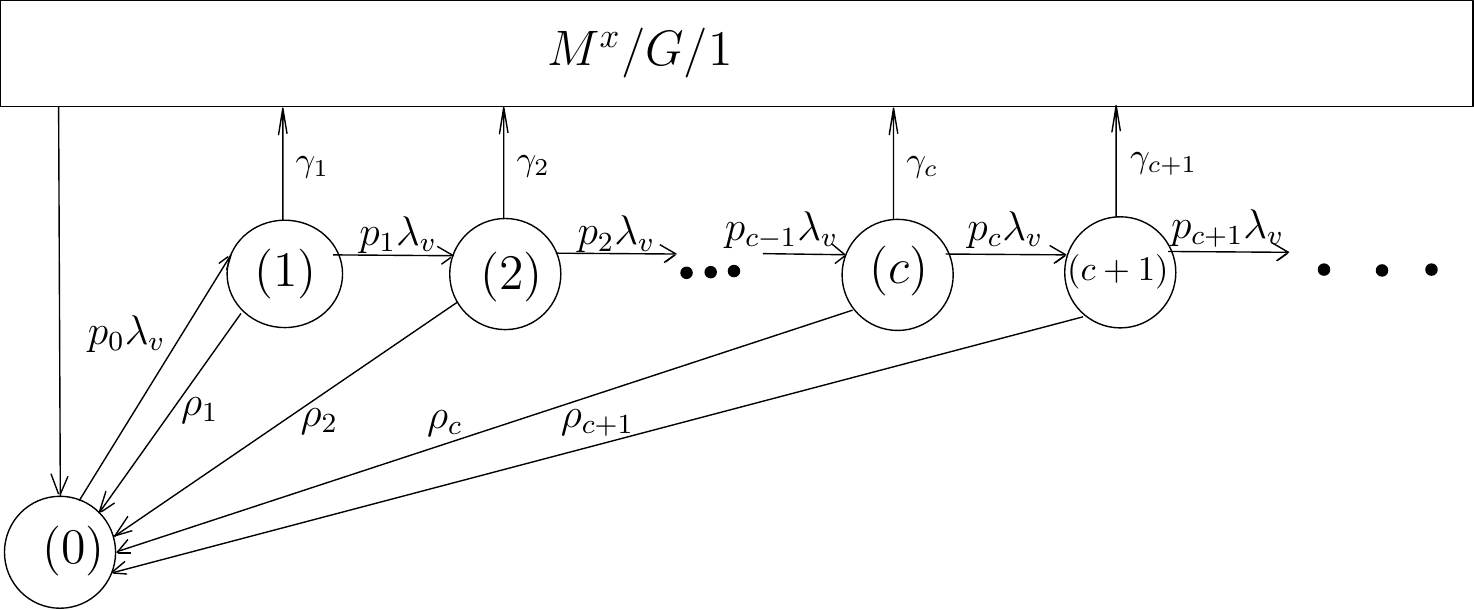}
	\caption{Rate transition diagram}
	\label{fig0n}
\end{figure}

 During the vacation mode
customers arrive in batches according to a  Poisson process with
 rate $\lambda_v$.  Let $g_{i}$ be the probability that the batch size is $i$, $i=1,...\,$.  A batch arrival that  sees $i$ customers admits  with probability $p_{i}$. During the vacation mode disasters may occur. At  moments of  disasters the system is cleared   \cite{YechDis}. The disaster rates are state dependent. That is,  when there are $i$ customers in the system disasters occur at rate $\rho_{i}$, i.e. the probability of a disaster in interval $dt$ is $\rho_{i}dt$.
When the number of customers in the system is $i$ the vacation terminates at rate $\gamma_{i}$, i.e. the probability that vacation terminates during an interval $dt$ is $\gamma_{i}dt$, $i=1,2,...$, ($\gamma_{0}=0$) (see Figure \ref{fig0n} for  the case where the batches are of size 1).

 Throughout this section it is assumed that  
  $\sum^{\infty}_{i=1}\gamma_i p_{(0,i)}=\xi<\infty$. In this case 
  \begin{equation}\label{eq.Yr}
  P(Y=i)=\frac{\gamma_{i} p_{(0,i)}}{\xi}
  \end{equation}
 and 

 \begin{equation}\label{eq.psiz}
 \Psi(z)=\frac{1}{\xi}\sum^{\infty}_{i=1}\gamma_i p_{(0,i)}z^i.
 \end{equation}
 

 
 
\textbf{a) Batches of size 1, i.e. $g_{1}=1$}

When $g_1=1$, 
the distribution of $Y$
 can be  expressed as a product form.
The stationary probabilities at vacation mode are given by

\begin{equation}\label{pok}
p_{(0,k)}=p_{(0,0)} \prod_{j=1}^k \frac{\lambda_v p_{j-1}}{\lambda_v p_j+\gamma_j+\rho_j}, \,\,k\geq 1
\end{equation}
and by (\ref{eq.Yr})
\begin{equation}\label{eq.PY}
P(Y=k)=\frac{1}{\xi}\gamma_{k} p_{(0,k)}, \,\,k\geq 1.
\end{equation}
Note that by (\ref{pok}) and (\ref{eq.PY}) the distribution of $Y$ {\color{red} is not a function of $p_{(0,0)}$} 

In this case we can obtain explicitly the steady state distribution of the number of customers in the vacation mode. 
%
To obtain $p_{(0,0)}$ we first find $p_{(0.)}$. Let $\mathcal{B}_{0}$ be the duration of the vacation mode in one cycle. The   equation for $E[\mathcal{B}_{0}]$ is given by

\begin{eqnarray}\label{EBo}
E[\mathcal{B}_{0}]
=&&\frac{1}{\lambda_{v}p_{0}}+\frac{1}{\lambda_{v}p_{1}+\gamma_{1}+\rho_{1}}+\sum_{i=2}^{\infty}\frac{1}{\lambda_{v}p_{i}+\gamma_{i}+\rho_{i}}\prod_{j=1}^{i-1}\frac{\lambda_{v}p_{j}}{\lambda_{v}p_{j}+\gamma_{j}+\rho_{j}}\nonumber\\
&& +	E[\mathcal{B}_{0}]\left(\frac{\rho_{1}}{\lambda_{v}p_{1}+\gamma_{1}+\rho_{1}}+
\sum_{i=2}^{\infty}\frac{\rho_{i}}{\lambda_{v}p_{i}+\gamma_{i}+\rho_{i}}\prod_{j=1}^{i-1}\frac{\lambda_{v}p_{j}}{\lambda_{v}p_{j}+\gamma_{i}+\rho_{j}}\right).
\end{eqnarray}
The 1st term on the right hand side  is the
 expected time in state $(0,0)$ and the 2nd  is the expected time in $(0,1)$. For the 3rd term, with probability
 $\prod_{j=1}^{i-1}\frac{\lambda_{v}p_{j}}{\lambda_{v}p_{j}+\gamma_{j}+\rho_{j}}$ state $(i,0)$ is reached before a disaster and before a transition to working mode. Next, the expected time in $(i,0)$ is $\frac{1}{\lambda_{v}p_{i}+\gamma_{i}+\rho_{i}}$. 
The term $\frac{\rho_{1}}{\lambda_{v}p_{1}+\gamma_{1}+\rho_{1}}$
 is the  probability  
 that a disaster occurs before a transition to the working mode or  to state $(0,2)$ starting at  state  $(0,1)$. Similarly, $\frac{\rho_{i}}{\lambda_{v}p_{i}+\gamma_{i}+\rho_{i}}\prod_{j=1}^{i-1}\frac{\lambda_{v}p_{j}}{\lambda_{v}p_{j}+\gamma_{j}+\rho_{j}}$ is the probability that the system reaches state $(0,i)$ before a disaster and before a  transition to the working mode and then a disaster occurs. Thus the expression in the brackets in the  second  line of (\ref{EBo}) is the probability of disaster before transition to working mode. Once a disaster occurs the time duration  in the vacation mode is distributed as $\mathcal{B}_{0}$. Solving for (\ref{EBo}) we obtain
	\begin{equation}
		E[\mathcal{B}_{0}]=\frac{\frac{1}{\lambda_{v}p_{0}}+\frac{1}{\lambda_{v}p_{1}+\gamma_{1}+\rho_{1}}+\sum_{i=2}^{\infty}\frac{1}{\lambda_{v}p_{i}+\gamma_{i}+\rho_{i}}\prod_{j=1}^{i-1}\frac{\lambda_{v}p_{j}}{\lambda_{v}p_{j}+\gamma_{j}+\rho_{j}}}{1-\left(\frac{\rho_{1}}{\lambda_{v}p_{1}+\gamma_{1}+\rho_{1}}+\sum_{i=2}^{\infty}\frac{\rho_{i}}{\lambda_{v}p_{i}+\gamma_{i}+\rho_{i}}\prod_{j=1}^{i-1}\frac{\lambda_{v}p_{j}}{\lambda_{v}p_{j}+\gamma_{j}+\rho_{j}}\right)}.
	\end{equation}\label{expectedvacation}
Obviously, the expected length of a cycle is:
	\begin{equation}
	E[T]=E[\mathcal{B}_{0}]+E[\mathcal{B}_{1}],
	\end{equation}
	where $\mathcal{B}_{1}$ is the duration of the working mode, where
	\begin{equation}
	E[\mathcal{B}_{1}]=\frac{E[S]E[Y]}{1-\lambda E[S]E[B]}.
	\end{equation} 
On the one hand
\begin{equation}
p_{0.}=\frac{E[\mathcal{B}_{0}]}{E[\mathcal{B}_{0}]+E[\mathcal{B}_{1}]}.
\end{equation}
On the other hand, from (\ref{pok})
\begin{equation}
p_{0.}=p_{(0,0)}\left(1+\sum_{k=1}^{\infty}\prod_{j=1}^k \frac{\lambda_{v} p_{j-1}}{\lambda_v p_j+\gamma_j+\rho_j}\right).
\end{equation}
Thus 
\begin{equation}
p_{(0,0)}\left(1+\sum_{k=1}^{\infty}
\prod_{j=1}^k \frac{\lambda_{v} p_{j-1}}{\lambda_v p_j+\gamma_j+\rho_j}\right)=\frac{E[\mathcal{B}_{0}]}{E[\mathcal{B}_{0}]+E[\mathcal{B}_{1}]}.
\end{equation}
Thus, we  find $p_{(0,0)}$ and by  (\ref{pok}) we obtain $p_{(0,k)}, k=1,2,...$ .

\textbf{b) A Markovian vacation with rational  parametric functions }

Here, 
we obtain
 (\ref{eq.psiz}) for a general batch size in the vacation mode. Recall that $g_{i}$ is the probability that the batch size is $i$, $i=1,\cdots,m$.
For general sequences $p_i$, $\gamma_i$, $\rho_i$ and $g_i$ it is too intricate   to find   $\Psi(z)$, since  (\ref{eq.psiz}) is the generating function of a product of two series.  Therefore

\begin{assumption}\label{assumptio2n}
	$p_j$, $\gamma_j$, $\rho_j$ and $g_j$ are  rational functions in  $j$, which mean that each of them can be expressed as a ratio of two polynomials in $j$.
\end{assumption}
The balance equations for $p_{(0,j)}$ are given by:
\begin{equation}\label{eq.po}
(\gamma_{j}+\lambda_v p_{j} +
\rho_{j})p_{(0,j)}=\sum_{i=\max(j-m,0)}^{j-1} \lambda_v p_{i} g_{j-i} p_{(0,i)},\,\,j=1,2,...
\end{equation}.

Then it can be shown by induction that $p_{(0,j)}, j \geq 1$ are also rational functions.  Equations (\ref{eq.po}) represent a special case of difference equations whose 
 solution is expressed in terms of hypergeometric functions as described in Knuth et al. \cite{Knuth}. 
The ordinary hypergeometric function is defined by
\begin{equation}\label{eq.F}
{}_k F_n(a_1,...,a_k;b_1,...b_n;z):=\sum_{i=0}^{\infty} \frac{(a_1)_i...(a_k)_i}{(b_1)_i...(b_n)_i}\frac{z^i}{i!},
\end{equation}
where $(a)_i=a \cdot (a+1) \cdot... \cdot (a+i-1)$. The hypergeometric  functions include many  special and elementary functions and usually appear as a solution of differential equations of special forms, (see e.g.  Petkovsek et al. \cite{Pet},  p.33).  

According to Knuth \cite{Knuth}, when the ratio $\frac{P(Y=j+1)}{P(Y=j)}$ is a rational function in j, then 
\begin{equation}
\Psi(z)=c\ _{k}F_{n}(a_1,a_2,...,a_k;b_1,b_2,...,b_n;z)z
\end{equation}
where $a_1,...,a_k,b_1,...,b_n,c$ are constants. $c=\frac{1}{F(a_1,a_2,...,a_k;b_1,b_2,...,b_h;1)}$ is the normalizing constant. In order to find the other constants we need to express the ratio $\frac{P(Y=j+1)}{P(Y=j)}$ as 

\begin{equation*}
\frac{P(Y=j+1)}{P(Y=j)}=\frac{(j+a_1)(j+a_2)...(j+a_s)}{(j+1)(j+b_1)(j+b_2)...(j+b_h)}.
\end{equation*}

As an example, let
 $\gamma_k={k+1}$, $p_k=\frac{1}{k+1}$, $\lambda=1$ and $g_1=1$. Then
\begin{equation}
\frac{P(Y=k+1)}{P(Y=k)}=\frac{\gamma_{k+1} p_{0,k+1}}{\gamma_k p_{0,k}}=\frac{(k+2)^2}{(k+1)^2 (k^2+4k+5)},
\end{equation}
and
\begin{equation}
\Psi(z)=\frac{z F(2,2;1,-2-i,-2+i;z)}{F(2,2;1,-2-i,-2+i;1)}.
\end{equation}
%
\subsubsection{Binomial reneging during the vacation }\label{sec.synchronised}
Theorem \ref{theorem.main} enables to generalize the model introduced by  Adan et al. \cite{IvoNew} from $M/G/1$ to  $M^{X}/G/1$  during the working mode. 
 The vacation period is the same as that    described in Section  \ref{sec.multiplevacation}, where $V$ is exponentially distributed($\gamma$).  During the vacation mode  customers arrive according to a Poisson process with rate $\lambda_{v}$ and there is no service. However,   
  customers may  abandon.  Abandonments occur according to  a  Poisson process($\xi$) independent of arrival process. At the moment of  abandonment  every  waiting customer  (if any) continues to wait   with probability $q=1-p$ or reneges with probability $p$, independently of the others. Thus the number of abandonments is a Binomial random variable out of the number of waiting customers, see figure \ref{MAE1} below.  
		\begin{figure}[h!]
			\centering
			\includegraphics[width=1\textwidth]{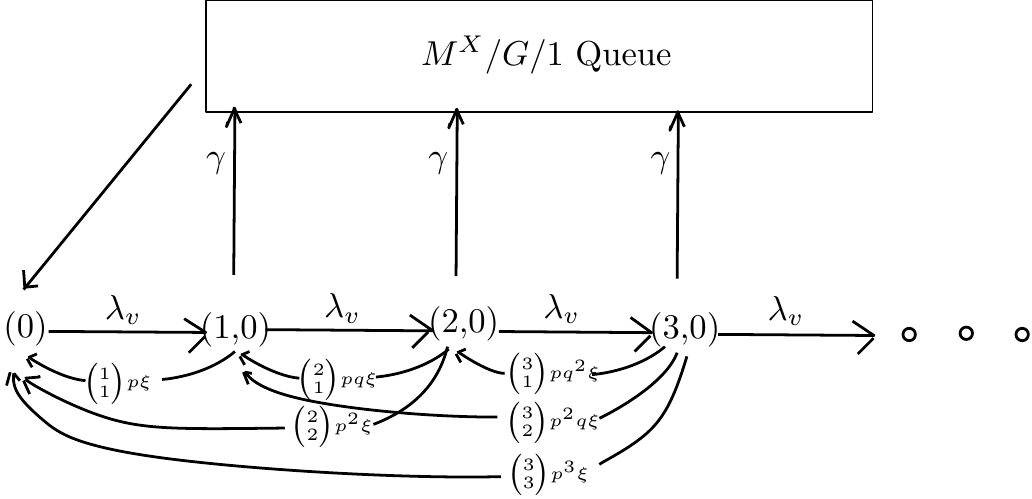}
				\caption{}
			\label{MAE1}
		\end{figure}
	 Adan et al. \cite{IvoNew}
 obtain the following  expression for the  partial  PGF of the number of customers in the vacation mode -- $G_0(z)$, where in the working mode the system runs as an M/M/1 queue   
 \begin{equation}\label{eq.G0synch1}
 G_0(z)=c \sum_{j=0}^{\infty}\prod_{k=0}^j\frac{\xi}{\gamma+\xi+\lambda (1-p)^k(1-z)},
 \end{equation}
 where  $c$ is the  normalizing  constant. 
 Clearly, the conditional PGF of the number of customers in Adan et al. \cite{IvoNew}  given that  the system runs in the vacation mode  is  the same as the PGF of  the number of customers here given that the system  runs in the vacation mode.  
 The PGF of the number of customers transferred from the vacation mode to the working  mode is 
 \begin{equation}\label{psizs}
 \Psi(z)=\frac{G_0(z)-G_0(0)}{G_0(1)-G_0(0)},
 \end{equation}
 where $ G_0(z)$ is given in (\ref{eq.G0synch1}), $p_{(0,0)}=G_{0}(0)$ and $p_{0.}=G_{0}(1)$. Obviously, $\Psi(z)$ is independent of the normalized constant $c$, since it is a proper PGF. 
The conditional  PGF of the number of customers given the working  is given by (\ref{decomposition}) with $\psi(z)$, which is the same as in (\ref{psizs}).

For the special case of the $M/M/1$ queue in the working mode we obtain

\begin{equation}
   \widetilde{G}_1(z) =\frac{(\mu-\lambda)z}{\mu-\lambda z}\frac{1-\Psi(z)}{(1-z)\Psi'(1)},
\end{equation}
where
\begin{equation}
\Psi'(1)=EY=\frac{G'_0(1)}{p_{0.}-p_{(0,0)}}=\frac{(1-p_{0.})(\mu-\lambda)}{\gamma (p_{0.}-p_{(0,0)})},
\end{equation}
which is the same as  in Adan et al. \cite{IvoNew}.  

\subsection{Related   systems with disasters}

Disasters are introduced in the literature with the convention that all the customers are cleared in one shot. That is,   when disasters occur the number of customers drops to zero. 
 The analysis of the applications introduced below is based  on a coupling argument.
Let $S_{I}$ be a queueing system with vacations as described in section \ref{secdescription} and assume that during the vacation the system runs as a certain  queueing system so-called  $S_{V}$. In the sequel we focus on special examples  of $S_{V}$. The transition rate $\gamma$ from  $S_{V}$ to the working mode is independent of the number of customers.  Next, let $S_{II}$ be an auxiliary queueing system without vacations that runs as a modified system $S_{V}$ with the added feature of disasters.  The disaster rate is independent of the number of customers and is equal to $\gamma$. We use the following 
\begin{criterion}\label{criterion}
In the  $S_{I}$ system, the conditional PGF of the number of customers given that the system runs in the vacation mode is equal to the  conditional PGF of the number of customers given that the server is busy in the  $S_{II}$ system. 
\end{criterion}
\begin{remark}
To  better understand the coupling idea    the following example is introduced.  	Suppose that  $S_{V}$ is a system with state dependent  disasters  that occur according to a Poisson process with rate $\xi_{i}$. Then,  in the system $S_{II}$ disasters occur according to a Poisson process with rate   $(\gamma+\xi_{i})$.
\end{remark}
 %
%
Let  $\pi_{i}$  be the  steady state probability of $i$ customers in $S_{II}$ and let $\Pi(z)$ be its PGF.
Recall that in $S_{I}$,   $Y$ is the number of customers transferred from the vacation mode to the working mode and $\Psi(z)$ is its PGF.
By Criterion \ref{criterion} 
\begin{equation} \label{Grel}
\Psi(z)=\frac{G_0(z)-p_{(0,0)}}{p_{0.}-p_{0,0}}=\frac{\Pi(z)-\pi_{0}}{1-\pi_{0}}.
\end{equation}
Once  $\Pi(z)$ is known, $\Psi(z)$ is also known.  Thus, by Theorem \ref{theorem.main} the conditional PGF of the number of customers given the working mode is found.

In the  special cases below  we apply (\ref{Grel}) and  Theorem \ref{theorem.main} to obtain the conditional PGF of  the  number of customers given the working mode in $S_{I}$ for special model versions  of $S_{V}$.

{\bf 	Special cases}
\begin{description}
 \item [1. M/M/1 queueing system  in vacation]
In the vacation mode the system behaves  as an $M/M/1$ queueing system with arrival rate $\lambda_{v}$ and service rate $\mu_{v}$ and the transition rate to working mode is  $\gamma$. For the  $M/M/1$ queue with disasters with rate $\gamma$ it is shown in  Kumar et al. \cite{Kumar} (Formula $3.2$ ) that
\begin{equation}\label{eq.pidis}
\pi_i=(1-\rho)\rho^i, i=0,1,...,
\end{equation}
where 
\begin{equation}
\rho=\frac{\gamma+\lambda_{v}+\mu_{v}-\sqrt{(\gamma+\lambda_{v}+\mu_{v})^2-4\lambda_{v} \mu_{v}}}{2 \mu_{v}}.
\end{equation}

 (\ref{Grel}) and (\ref{eq.pidis}) imply that  
\begin{equation}
\Psi(z)=\frac{(1-\rho)z}{1-\rho z}.
\end{equation}


\item[2. Birth and death process (BDP) with chain sequence rates]
Consider a queueing system where during vacation the system's behavior is that of BDP with chain sequence rates. 
 Lenin et al.  \cite{Lenin} 
 define  a birth-death  system whose birth and death rates are given by the following chain sequence:
$\lambda_0=1$,  $\lambda_n+\mu_n=1$, $\lambda_{n-1} \mu_n=a$  for $n=1,2,...$, where   $0<a\leq 0.25$ is a given constant.
For example, the case  $\lambda_n=\mu_n=\frac{1}{2}$ is a BDP system with $a=0.25$. In general, for $n=1,2,3,...$
\begin{equation}
\begin{array}{lll}
	\lambda_n=\frac{\alpha U_{n+1}(\frac{1}{\alpha})}{2 U_{n}(\frac{1}{\alpha})},&&
	\mu_{n}=\frac{\alpha U_{n-1}(\frac{1}{\alpha})}{2 U_{n}(\frac{1}{\alpha})}.
\end{array}
\end{equation}

Here, $\lambda_0=1$, $\mu_0=0$, $U_n(\cdot)$ is the second type  Chebyshev's  polynomial of order $n$ (see  \cite{Mason}) and $\alpha=2 \sqrt{a}$.
 $U_n(\cdot)$ is defined recursivly by
\begin{equation}
    U_n(x)=\frac{sin((n+1)\theta)}{sin(\theta)},n=0,1,2,3,...;cos(\theta)=x.
\end{equation}
For example, the case $a=0.2$ yields
$\lambda_1=0.796,\lambda_2=0.746,\lambda_3=0.729,\lambda_4=0.723,\lambda_5=0.721,...$ ,  so $\lambda_n\geq \mu_n$ and $\lambda_n$ decreases  and $\mu_n$ increases. 
This system is not stable, since  $\lambda_n\geq\mu_n$.
However, if disasters occurring at a fixed rate $\gamma$, are added to the model  the system is  stable and its steady state probabilities are given in  Kumar et al. (\cite{Kumar}) 
\begin{equation}\label{piQBD}
\pi_n=\frac{\gamma}{2\sqrt{a}}U_n\left(\frac{1}{2\sqrt{a}}\right)\left(\frac{\gamma+1-\sqrt{(\gamma+1)^2-4a}}{2\sqrt{a}}\right)^{n+1},n=0,1,2....
\end{equation}
$\Psi(z)$ is obtained by  substituting (\ref{piQBD}) in  (\ref{Grel}).





%


\item[3. Binomial reneging and disasters]

%
Consider a system that runs as  described in Subsection \ref{sec.synchronised} with the added feature of disasters occurring according to a Poisson process with rate $\gamma$.   
\begin{figure}[h!]
	\centering
	\includegraphics[width=1\textwidth]{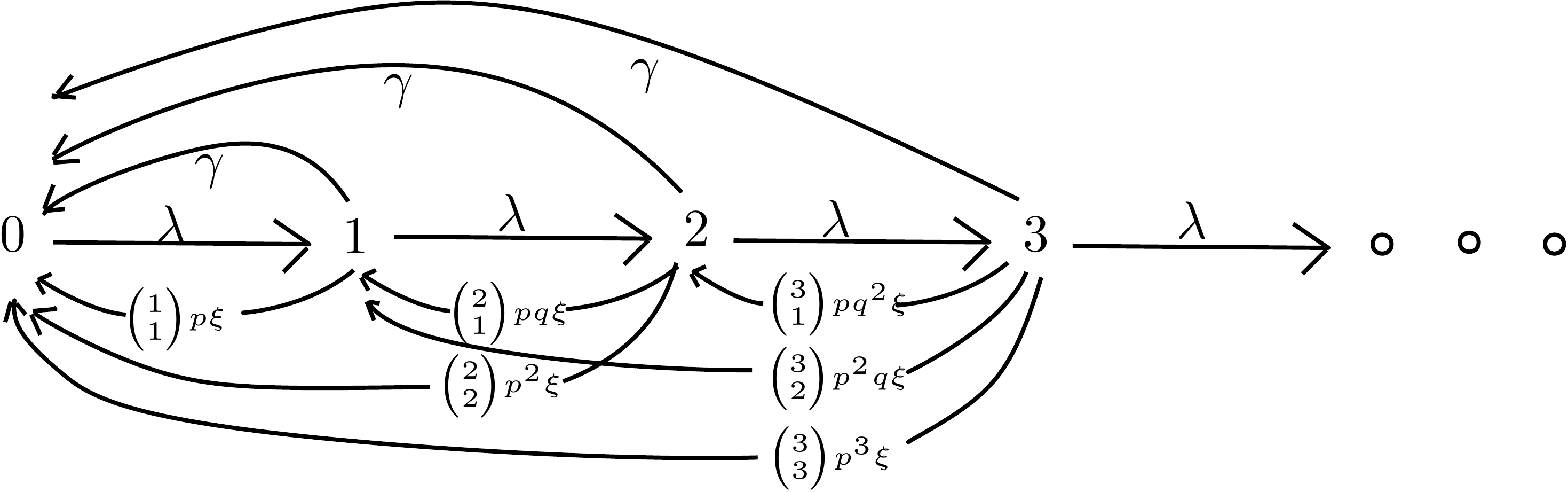}
	\caption{Synchronized services with disasters - rate transition diagram}
	\label{fig16}
\end{figure}
%
 The transition diagram of this model is presented in Figure \ref{fig16} below. 

By  (\ref{Grel}) 

%
%

\begin{equation}
    \Pi(z)=\frac{G_0(z)}{p_{0.}}=\frac{G_0(z)}{G_0(1)} 
\end{equation}
where, $G_0(z)$ is given in (\ref{eq.G0synch1}) .

\item[4.  $\mathbf{ M^{X}/G/1}$ with disasters] Mytalas and Zazanis \cite{Mytalas} introduce  an $M^{X}/G/1$ queue with disasters.  When state 0 is reached  by a  service completion the server takes a vacation and when state  0 is reached by a disaster a repair starts. Customers continue to arrive during the vacation or the repair periods, but  the  service is stopped during these periods.

In our model we  consider an $M^{X}/G/1 $ queue in the working mode. During the vacation mode the system runs as   another $M^{X}/G/1$ queue with Poisson arrivals at rate $\lambda_{V}$, generic service time $S_{V}$ with LST $\widetilde{F}_{S_{V}}$ and the PGF of the batch is $B_{V}(z)$.   During the vacation disasters may occur at rate $\xi$. Recall that the transition rate from the vacation mode to the working mode is $\gamma$.
By (\ref{Grel}) $\Psi(z)$ is 
 the conditional PGF of the number of customers in the queueing system introduced by \cite{Mytalas} with disaster rate $\gamma+\xi$   given that the server is busy.  Thus the PGF is
\begin{equation}\label{Psizzaz}
	\Psi(z)=\frac{\xi+\gamma}{1-B_{V}(z_{\xi+\gamma})}\frac{z(B_{V}(z_{\xi+\gamma})-B_{V}(z))}{\widetilde{F}_{S_{V}}(\xi+\gamma+\alpha(z))-z}\frac{1-\widetilde{F}_{S_{V}}(\xi+\gamma+\alpha(z))}{\xi+\gamma+\alpha(z)},
\end{equation}
where $z_{\xi +\gamma }$ is the unique solution of the equation 
\[z=\widetilde{F}_{S_{V}}(\xi+\gamma+\lambda_{V}-\lambda_{V}B_{V}(z))\]
and
\[\alpha(z)=\lambda_{V}(1-B_{V}(z)).\]
\end{description}
By substituting (\ref{Psizzaz})     in (\ref{decomposition}) we obtain the conditional PGF of the number of  customers given the working mode.

In addition,  the steady-state probability that the system is in the vacation mode is given by
\[p_{0.}=\frac{\Phi_{V}'(0)}{\Upsilon'(0)+\Phi'(0)},
\]
where $\Phi_{V}(\alpha)$ and $\Upsilon(\alpha)$ are the LSTs of the vacation period and the working  period, respectively, see Kleiner et al. \cite{kleiner}. 
Multiplying the conditional PGF of the number of customers given the working mode by $1-p_{0.}$ we obtain the unconditional PGF of the number of customers in the working mode. 
Note that $\Psi(z)$ is also the PGF of the number of customers in  steady state in the vacation mode  given that the server is  busy. Thus $p_{0,0}+(p_{0.}-p_{0,0})\Psi(z))$
  is the PGF of the number of customers in steady state in vacation, where $p_{0,0}$ is also  given in Kleiner et al. \cite{kleiner}.

		%
		
		\end{document}